\documentclass[letterpaper]{sig-alternate-per-wg73}
\usepackage{amssymb,mathrsfs,verbatim,color,lscape,xspace,enumerate,enumitem}
\allowdisplaybreaks
\usepackage{listings,mathtools}
\newtheorem{theorem}{Theorem}[section]
\newdef{remark}{Remark}
\newtheorem{corollary}{Corollary}
\newtheorem{definition}{Definition}

\newcommand{\splus}{{+}}  
\newcommand{\sgreater}{{>}}
\newcommand{\ssmaller}{{<}}

\newcommand{\pr}[1][]{\ensuremath{\mathbb{P}{#1}}\xspace}   
\newcommand{\apr}[1]{\ensuremath{\hat{#1}}}
\newcommand{\fun}[2][s]{\ensuremath{{#2}(#1)}}     
\newcommand{\diag}{\ensuremath{\text{diag}}\xspace}  
\newcommand{\adj}{\ensuremath{\text{adj}}\xspace}  
\newcommand{\rank}{\ensuremath{\text{rank}}\xspace}  
\newcommand{\mean}[1][]{\ensuremath{\mu_{#1}}\xspace}    

\newcommand{\ptsd}[1][t]{\ensuremath{F_{p}(#1)}\xspace}        
\newcommand{\htsd}[1][t]{\ensuremath{F_{h}(#1)}\xspace}        
\newcommand{\gsd}[1][t]{\ensuremath{F(#1)}\xspace}        
\newcommand{\gsdc}[2][t]{\ensuremath{G_{#2}(#1)}\xspace}        

\newcommand{\ltpts}[1][]{\ensuremath{\widetilde{F}_{p}^{#1}(s)}\xspace}        
\newcommand{\lthts}[1][]{\ensuremath{\widetilde{F}_{h}^{#1}(s)}\xspace}        
\newcommand{\ltgs}[1][]{\ensuremath{\widetilde{F}^{#1}(s)}\xspace}        
\newcommand{\ltw}[1][s]{\ensuremath{\widetilde{v}(#1)}\xspace}        
\newcommand{\ltwmix}[1][s]{\ensuremath{\widetilde{v}_\epsilon(#1)}\xspace}        
\newcommand{\vltw}[1][s]{\ensuremath{\widetilde{\mathbf{\Phi}}(#1)}\xspace}        
\newcommand{\ltwc}[2][]{\ensuremath{\widetilde{\phi}_{#2}^{#1}(s)}\xspace}        

\newcommand{\mxtrans}{\ensuremath{\mathbf{P}}\xspace}           
\newcommand{\mxrates}{\ensuremath{\mathbf{\Lambda}}\xspace}           
\newcommand{\mxprobs}{\ensuremath{\mathbf{Q}}\xspace}           
\newcommand{\ltmgs}[1][s]{\ensuremath{\widetilde{\mathbf{G}}(#1)}\xspace}     
\newcommand{\ltmgsmix}[1][s]{\ensuremath{\widetilde{\mathbf{G}}_\epsilon(#1)}\xspace}     
\newcommand{\uv}{\ensuremath{\mathbf{e}}\xspace}           
\newcommand{\im}[1][]{\ensuremath{\mathbf{I}_{#1}}\xspace}           
\newcommand{\vup}{\ensuremath{\mathbf{u}}\xspace}           
\newcommand{\vupmix}{\ensuremath{\mathbf{u}_\epsilon}\xspace}           
\newcommand{\va}[1][]{\ensuremath{\mathbf{a}_{#1}}\xspace}           

\newcommand{\rootp}[1]{\ensuremath{s_{#1}}\xspace}    
\newcommand{\rootpmix}[1]{\ensuremath{s_{\epsilon,#1}}\xspace}    
\newcommand{\rootnden}[1]{\ensuremath{x_{#1}}\xspace}    
\newcommand{\rootnnum}[1]{\ensuremath{y_{#1}}\xspace}    
\newcommand{\workload}[1][]{\ensuremath{V_{#1}}\xspace}    

\newcommand{\epts}[1][]{\ensuremath{B^e_{#1}}\xspace}                 
\newcommand{\ehts}[1][]{\ensuremath{C^e_{#1}}\xspace}                 

\begin{document}

\title{Corrected phase-type approximations for the workload of the MAP/G/1 queue with heavy-tailed service times}
\author{
    Eleni Vatamidou\\
      \begin{footnotesize}
        EURANDOM and
      \end{footnotesize}\\
    \begin{footnotesize}
        Eindhoven University
    \end{footnotesize}\\
    \begin{footnotesize}
        of Technology
    \end{footnotesize}\\
    e.vatamidou@tue.nl\\
    \and
    Ivo Adan\\
          \begin{footnotesize}
        EURANDOM and
      \end{footnotesize}\\
    \begin{footnotesize}
        Eindhoven University
    \end{footnotesize}\\
    \begin{footnotesize}
        of Technology
    \end{footnotesize}\\
    i.j.b.f.adan@tue.nl\\
    \and
    Maria Vlasiou\\
      \begin{footnotesize}
        EURANDOM, CWI, and
      \end{footnotesize}\\
    \begin{footnotesize}
        Eindhoven University
    \end{footnotesize}\\
    \begin{footnotesize}
        of Technology
    \end{footnotesize}\\
    m.vlasiou@tue.nl\\
    \and
    Bert Zwart\\
          \begin{footnotesize}
        EURANDOM, CWI,
      \end{footnotesize}\\
      \begin{footnotesize}
        VU University Amsterdam,
      \end{footnotesize}\\
      \begin{footnotesize}
      Georgia Institute of
      \end{footnotesize}\\
          \begin{footnotesize}
       Technology, and Eindhoven
    \end{footnotesize}\\
    \begin{footnotesize}
      University of Technology
    \end{footnotesize}\\
    Bert.Zwart@cwi.nl}

\maketitle

\begin{abstract}
In many applications, significant correlations between arrivals of load-generating events make the numerical evaluation of the load of a system a challenging problem. Here, we construct very accurate approximations of the workload distribution of the MAP/G/1 queue that capture the tail behavior of the exact workload distribution and provide a small relative error. Motivated by statistical analysis, we assume that the service times are a mixture of a phase-type and a heavy-tailed distribution. With the aid of perturbation analysis, we derive our approximations as a sum of the workload distribution of the MAP/PH/1 queue and a heavy-tailed component that depends on the perturbation parameter. We refer to our approximations as {\it corrected phase-type approximations}, and we exhibit their performance with a numerical study.
\end{abstract}

\begin{center}
\end{center}

\section{Introduction}\label{intro}
The evaluation of the workload of a MAP/G/1 queue is an important problem that has been widely studied in the literature. For an extensive review see \cite{lucantoni93-BMAP}. Although closed-form expressions for the evaluation of the workload are available, they are practical only in the case of phase-type (PH) service times. When the workload distribution cannot be computed exactly, it needs to be approximated. Here, we develop a new method to construct reliable approximations for the workload distribution for heavy-tailed service times.

Two main directions for approximating the workload distribution are PH and asymptotic approximations. When the service times are light-tailed, a common approach to approximate the workload with high accuracy is by approximating the service time distribution with a PH one \cite{feldmann98,starobinski00}. We refer to these methods as {\it phase-type approximations}, because the approximate workload distribution has a PH representation \cite{ramaswami90}. However, in many cases, a heavy-tailed distribution is most appropriate to model the service times \cite{embrechts-MEE,rolski-SPIF}. For the special class of subexponential service times, asymptotic approximations for the workload distribution are available, which provide a good fit only at the tail \cite{asmussen-RP,embrechts82}.

In this paper, we develop approximations of the workload distribution for heavy-tailed service times that maintain the computational tractability of PH approximations, capture the correct tail behavior and provide small absolute and relative errors. Also, they have the advantage that finite higher-order moments for the service times are not required. In order to achieve these desirable characteristics, our key idea is to use a mixture model for the service times.

The idea of our approach stems from fitting procedures of the service time distribution to data. Heavy-tailed statistical analysis suggests that only a small fraction of the upper-order statistics of a sample is relevant for estimating tail probabilities \cite{resnick-HTP}. The remaining data set may be used to fit the bulk of the distribution. Since PH distributions are dense in the class of all positive definite probability distributions \cite{asmussen-APQ}, a natural choice is to fit a PH distribution to the remaining data set \cite{asmussen96}. As a result, a mixture model for the service times is a natural assumption.

In short, we consider the service time distribution as a mixture of a PH distribution and a heavy-tailed one. As ``base" model we use the model appearing when all heavy-tailed customers are removed, and we interpret the heavy-tailed term of the mixture model as perturbation of the PH one. Using perturbation analysis, we find our approximations for the workload in the mixture model as a sum of the workload of the base model and a heavy-tailed component that depends on the perturbation parameter. In a previous study (cf. \cite{vatamidou13}), we carried out this project for Poisson arrivals. Here we develop an extension to MAP's.

The rest of the paper is structured as follows. In Section~\ref{s.model}, we introduce the notation for the model under consideration, and in Section~\ref{s.algorithm}, we present the algorithm to construct our approximations. In Section~\ref{s.approximations}, we give their formulas and we also specialize to the M/G/1 queue. Finally, in Section~\ref{s.numerical example}, we perform an illustrative numerical experiment.

\section{Model description}\label{s.model}
We consider a single server queue with FIFO discipline, where customers arrive according to a Markovian Arrival Process (MAP) with $N$ states. The regulating Markov chain $\{Z_n\}_{n\geq 0}$ has an irreducible transition probability matrix \mxtrans and stationary distribution $\boldsymbol\pi$. Transitions from state $i$ occur at an exponential rate $\lambda_i$ (directed to state $j$ with probability $p_{ij}$). With probability $q_i$, a transition from $i$ corresponds to an arrival of a (real) customer with service time distribution $\gsd$ (independent of the state), and otherwise, it corresponds to an arrival of a (dummy) customer with zero service time. So the service time distribution of a customer arriving in state $i$ is $\gsdc{i} = q_i\gsd + 1-q_i$. Observe that if the service time distribution of a real customers \gsd was depending on the state of the system then we would allow for cross-correlations between the arrivals and the services, which is not the case in our model.

In matrix form, the above quantities can be written as $\mxprobs = \diag(q_1,\dots,q_N)$, $\mxrates = \diag(\lambda_1,\dots,\lambda_N)$ and $\ltmgs = \ltgs \mxprobs + (\im - \mxprobs)$, where \ltgs denotes the Laplace-Stieltjes transform (LST) of the service time distribution \gsd of a real customer, and \im stands for the identity matrix with dimension $N$. Finally, let \mean be the mean of the service time distribution \gsd. We assume that the system is stable, namely $\boldsymbol\pi \left(\mxrates^{-1} - \mean \mxprobs\right) \uv>0$, where \uv is the column vector with all elements equal to 1. If \ltwc{i} denotes the LST of the steady-state workload in state $i$, the following theorem holds for the transform vector $\vltw = [\ltwc{1},\dots,\ltwc{N}]$ (cf. Th. 3.1 in \cite{adan03}).

\begin{theorem}\label{t3.equations for the transform vector}
  If the system is stable, there exists a unique vector $\vup=\left[u_1,...,u_N\right]$, such that \vltw satisfies
  \begin{align}
    \vltw \left(\ltmgs \mxtrans \mxrates +s \im -\mxrates \right) &= s \vup, \notag \\ 
    \vltw[0] \uv &= 1. \notag 
  \end{align}
\end{theorem}

To determine the unknown vector \vup we have (cf. Th. 3.2 \& 3.3. in \cite{adan03}):
\begin{theorem}\label{t3.solution for vector u}
  It holds that
  \begin{enumerate}
   \item $\det\left(\ltmgs \mxtrans \mxrates +s\im -\mxrates\right)=0$ has exactly $N$ roots $\rootp{i}$, with $\rootp{1}=0$ and $Re(\rootp{i})>0$, $i=2,\dots,N$.\label{part1}
    \item Let $\va[i]$ be non-zero column vectors satisfying
        \begin{equation*}
          \left(\ltmgs[\rootp{i}]\mxtrans \mxrates +\rootp{i} \im -\mxrates \right)\va[i] = 0, \quad i=2,\dots,N.
        \end{equation*}
        Then, provided all \rootp{i} are distinct, \vup is the unique solution to the $N$ linear equations:
        \begin{align*}
          \vup \mxrates^{-1} \uv &= \boldsymbol\pi\left(\mxrates^{-1} - \mean \mxprobs\right) \uv, \\
          \vup \va[i] &= 0, \qquad \qquad i=2,\dots,N.
        \end{align*}
  \end{enumerate}
\end{theorem}

Combining the results of Theorems~\ref{t3.equations for the transform vector} and \ref{t3.solution for vector u}, the LST of the total workload $V$ in the system is
\begin{equation}\label{e3.Laplace transform workload as a fraction}
  \ltw = \frac{s \cdot \vup \cdot \adj \left(\ltmgs \mxtrans \mxrates +s \im -\mxrates\right) \uv}{\det \left(\ltmgs \mxtrans \mxrates +s \im -\mxrates\right)},
\end{equation}
where \adj denotes the adjoint of a square matrix.

We now assume that the service time distribution of a real customer has the form $\gsd = (1-\epsilon)\ptsd + \epsilon \htsd$, $\epsilon \in [0,1)$, for some PH (\ptsd) and heavy-tailed (\htsd) distributions, with finite means \mean[p] and \mean[h], respectively. Theorem~\ref{t3.solution for vector u} guarantees that the RHS of \eqref{e3.Laplace transform workload as a fraction} in this model with mixed service time distribution is well-defined in the positive half-plane. However, if the LST \lthts does not have a closed-form expression (e.g.Pareto), Laplace inversion of \eqref{e3.Laplace transform workload as a fraction} cannot be applied to find the distribution of the workload $V_\epsilon$ in this mixture model.

In the next section, we describe how to create an approximation for $V_\epsilon$, by approximating its LST.

\section{Approach}\label{s.algorithm}
The steps to construct our approximations are:
\begin{enumerate}
  \item Use a PH approximation as base model.
    \begin{enumerate}
      \item Set $\epsilon=0$ and $\ltpts = \fun{q_n}/\fun{p_m}$, where \fun{q_n} and \fun{p_m} are polynomials of degrees $m$ and $n$ respectively, with $n\leq m-1$, so that \ltpts is the LST of a PH-distribution.
      \item Use Theorem~\ref{t3.solution for vector u} to determine the vector \vup and find the adjoint matrix $\adj \left(\ltmgs \mxtrans \mxrates +s \im -\mxrates\right)$.
      \item Find the LST \ltw (see \eqref{e3.Laplace transform workload as a fraction}) as
        \begin{equation}\label{e3.L-T waiting final form}
            \ltw = \frac{\vup \uv \prod_{j=1}^{mr} (s+\rootnnum{j})}{\prod_{j=1}^{mr} (s+\rootnden{j})},
        \end{equation}
        where $Re(\rootnnum{j})>0$, $Re(\rootnden{j})>0$, $j=1,\dots,mr$, and $r$ is an non-negative integer smaller than or equal to the $\rank$ of \mxtrans.
      \item Apply Laplace inversion to \eqref{e3.L-T waiting final form} to find analytically $\pr(V>t)$. \label{step}
    \end{enumerate}
  \item Find the parameters of the mixture model as perturbation of the base model's ones (parameters affected by the perturbation bear an index $\epsilon$).
    \begin{enumerate}
      \item For the matrix $\ltmgsmix \mxtrans \mxrates +s \im -\mxrates$, find its determinant and its adjoint matrix .
      \item Evaluate the vector \vupmix and the roots \rootpmix{i}, $i=1,\dots,N$, using an extension of Theorem~\ref{t3.solution for vector u} (omitted due to space limitations).
    \end{enumerate}
  \item Find the LST of the workload $V_\epsilon$ as perturbation of \ltw, by keeping only up to $\epsilon$-order terms, i.e.,
    \begin{equation}\label{e3.laplacetransformperturbed}
      \ltwmix = \ltw + \epsilon \ltw \fun[s]{k} + O(\epsilon^2),
    \end{equation}
    where \fun[s]{k} is well-defined for positive values.
\end{enumerate}

Our proposed approximations are constructed by applying Laplace inversion to the up to $\epsilon$-order terms of \eqref{e3.laplacetransformperturbed}. In the next section, we give their formulas.

\section{Corrected PH approximations}\label{s.approximations}
Let \epts and \ehts be the generic stationary excess PH and heavy-tailed service times, respectively. Moreover, let $E_\lambda$ be an exponential r.v.\ with rate $\lambda$, and let $\workload'$ be independent and follow the same distribution of \workload.\ The next theorem shows that each term in the Laplace inverse $\mathcal{L}^{-1}\{\ltw \fun[s]{k}\}$ has a probabilistic interpretation.
\begin{theorem}\label{t3.perturbed workload distribution}
  There exist unique coefficients $\alpha$, $\beta$, $\gamma$, $\alpha_j$, $\beta_j$, $\gamma_j$, $j=1,\dots,mr$ and $\delta_i$, $\eta_i$, $\theta_i$, $i=2,\dots,N$, s.t.
  \begin{footnotesize}
  \begin{align*}
    &\mathcal{L}^{-1}\{\ltw \fun[s]{k}\} = \frac1{\vup \uv}\Bigg( \beta \big( \mean[p]\pr(\workload \splus \epts \sgreater t)
                     - \mean[h]\pr(\workload \splus \ehts \sgreater t) \big) \\
                    +& \sum_{j=1}^{mr}\beta_j \big( \mean[p] \pr(\workload \splus \epts \splus E_{y_j} \sgreater t) - \mean[h] \pr(\workload \splus \ehts \splus E_{y_j} \sgreater t) \big) \\
                    +& \sum_{i=2}^N \eta_i \big( \mean[p] \pr(t \ssmaller \workload \splus \epts \ssmaller t \splus E_{\rootp{i}})
                     - \mean[h] \pr(t \ssmaller \workload \splus \ehts \ssmaller t\splus E_{\rootp{i}}) \big)\\
                    +&\gamma \big( \mean[p]\pr(\workload \splus \workload' \splus \epts \sgreater t)
                     - \mean[h]\pr(\workload \splus \workload' \splus \ehts \sgreater t) \big) \\
                    +& \sum_{j=1}^{mr} \gamma_j \big( \mean[p] \pr(\workload \splus \workload' \splus \epts \splus  E_{y_j} \sgreater t)
                     - \mean[h] \pr(\workload \splus \workload' \splus  \ehts \splus E_{y_j} \sgreater t)\big)\\
                    +& \sum_{i=2}^N \theta_i \big( \mean[p] \pr(t \ssmaller \workload \splus \workload' \splus \epts \ssmaller t\splus E_{\rootp{i}}) \\
                     &- \mean[h] \pr(t \ssmaller \workload \splus \workload' \splus \ehts \ssmaller t\splus E_{\rootp{i}}) \big)
                     +\alpha \pr(\workload \sgreater t)\\
                    &+ \sum_{j=1}^{mr}\alpha_j \pr(\workload \splus E_{y_j}> t)
                     + \sum_{i=2}^N \delta_i \pr(t \ssmaller \workload \ssmaller t \splus E_{\rootp{i}})\Bigg).
  \end{align*}
  \end{footnotesize}
\end{theorem}

\begin{remark}
  In Theorem~\ref{t3.perturbed workload distribution}, we assumed for convenience that all \rootnnum{j} are simple, but the result can be generalized to roots with multiplicity greater than one. Also, the unique coefficients are found in a straightforward way, but we omit the details.
\end{remark}

\begin{remark}
  We assumed that all \rootp{i} and \rootnnum{j} are real-valued. If e.g.\ \rootp{2} is complex, then we write $E_{Re(\rootp{2})}$ instead of $E_{\rootp{2}}$. The imaginary part cancels out when we combine each complex root with its conjugate.
\end{remark}

Let $\apr{V}_\epsilon$ denote the approximation of $V_\epsilon$. Following the previous result, we have
\begin{definition}
  The {\it corrected} PH approximation is
  \begin{equation}
    \pr(\apr{V}_\epsilon>t) := \pr(\workload > t) + \epsilon \mathcal{L}^{-1}\{\ltw \fun[s]{k}\},
  \end{equation}
  where $\mathcal{L}^{-1}\{\ltw \fun[s]{k}\}$ is given by Theorem~\ref{t3.perturbed workload distribution} and $\pr(V>t)$ follows from step \eqref{step}.
\end{definition}

For the M/G/1, the coefficients of Theorem~\ref{t3.perturbed workload distribution} are directly obtained from the parameters of the system.
\begin{corollary}
  In case of Poisson arrivals with rate $\lambda$, the {\it corrected} PH approximation takes the form
  \begin{align*}
    \pr(\apr{V}_\epsilon &>t) := \pr(\workload > t) + \epsilon\frac{\lambda}{1-\lambda \mean[p]} \Big( (\mean[p]-\mean[h])\pr(\workload > t) \\
                             +& \mean[h]\pr(\workload \splus \workload' \splus \ehts \sgreater t) - \mean[p] \pr(\workload \splus \workload' \splus \epts \sgreater t) \Big). \label{e.corected replace}
  \end{align*}
\end{corollary}

Using a test model, in Section~\ref{s.numerical example}, we check the performance of our approximations.

\section{Numerical example}\label{s.numerical example}
We consider a MAP with Erlang-2 distributed interarrival times, where the exponential phases have both rate $\lambda$, namely $\lambda_i=\lambda$, $i=1,2$. For the service times, we use a mixture of an Exp$(\nu)$ distribution and a heavy-tailed one (cf.\ \cite{abate99a}) with LST $\lthts = 1-\frac{s}{(\kappa+\sqrt{s})(1+\sqrt{s})}$. The exact workload distribution for this mixture of distributions can be calculated by following a similar idea as the proof of Th.~9 in \cite{vatamidou13}. For our numerical examples, we select $\kappa=2$, $\nu=3$ and $\epsilon=0.01$.

\begin{table}[h]
\centering
    \begin{tabular}{|r|clc|}\hline
      $t$   &exact           & PH                 & corrected-PH \\ \hline
      0     & 0.837500       & 0.833333                   & 0.837213 \\
      5     & 0.061452       & 0.031781                   & 0.060882 \\
      10    & 0.023269       & 0.001212                   & 0.023544 \\
      15    & 0.017579       & 0.000046                   & 0.017862 \\
      20    & 0.014979       & $1.76 \times 10^{-6}$      & 0.014091 \\
      25    & 0.013301       & $6.72 \times 10^{-8}$      & 0.013126 \\
      30    & 0.012090       & $2.56 \times 10^{-9}$      & 0.011867 \\
      35    & 0.011162       & $9.78 \times 10^{-11}$     & 0.010943 \\
      40    & 0.010419       & $3.73 \times 10^{-12}$     & 0.010220 \\
      45    & 0.009809       & $1.42 \times 10^{-13}$     & 0.009601 \\
      50    & 0.009294       & $5.42 \times 10^{-15}$     & 0.009106 \\ \hline
    \end{tabular}
  \caption{Tail probabilities of the exact workload, the PH and the {\textit{corrected}} PH approximations for $\epsilon=0.01$ and load $0.8375$.}\label{table-error}
\end{table}

From Table~\ref{table-error}, we observe that the {\it corrected} PH approximation yields a significant improvement to its PH counterpart. The difference between the exact tail probabilities of the workload and the {\it corrected} PH approximation is $O(10^{-4})$, while for the PH approximation it is $O(10^{-2})$, for all values. The magnitude of the improvement we achieve with the {\it corrected} PH is evident by looking at the relative errors of the involved approximations. The relative error of the PH easily reaches values close to 1, while the {\it corrected} PH gives a relative error $O(\epsilon)$.

\end{document}